\documentclass[12pt,a4paper]{article}

\usepackage{amssymb}
\usepackage{amsmath}
\usepackage{amsthm}
\usepackage{xypic}
\theoremstyle{plain}
\newtheorem{theorem}{Theorem}
\newtheorem{proposition}[theorem]{Proposition}

\theoremstyle{remark}

\theoremstyle{definition}

\def\varinjlim_#1{\lim\limits_{\longrightarrow\atop{#1}}}

\def\id{\mathop{\rm id}\nolimits}

\def\diag{\mathop{\rm diag}\nolimits}
\def\pt{\mathop{\rm pt}\nolimits}
\def\Ob{\mathop{\rm Ob}\nolimits}
\def\Mor{\mathop{\rm Mor}\nolimits}
\def\I{\mathop{\rm I}\nolimits}

\begin{document}
\author{A.V. Ershov}
\title{Supplement 2\\ to the paper ''Floating bundles and their
applications''}
\date{}
\maketitle This paper is the supplement to the section 2 of the
paper "Floating bundles and their applications" \cite{pap1}. Below
we study some properties of category, connected with
cobordism rings of FBSP. In particular, we shall show that
it is the tensor category.

In \cite{pap1} the series ${\frak G}(x,y)\in H[[x,y]]=\Omega^*_U
(\widetilde{Gr}),$ where $H=\Omega_U^*(Gr),$ was defined.
Recall that it corresponds to the direct limit $\kappa$
of the maps $\kappa_{k,l}\colon \widetilde{Gr}_{k,kl}\rightarrow
\mathbb{C}P^{kl-1},$ where $\widetilde{Gr}_{k,kl}$
is the canonical FBSP over $Gr_{k,kl}$ ($(k,l)=1$).
In \cite{pap1} some properties of ${\frak G}(x,y)$ were
studied. In particular, it was shown that
$$
(\varepsilon{\frak G})(x,y)=F(x,y),
$$
where $\varepsilon \colon H\rightarrow R=\Omega^*_U(\pt)$
is the counit of the Hopf algebra $H$ and $F(x,y)\in R[[x,y]]$
is the formal group of geometric cobordisms.

Let $\varphi_{k,l}$ be the map
$$
\kappa_{k,l}\times \id_{\widetilde{Gr}_{k,kl}}\colon
\widetilde{Gr}_{k,kl}\rightarrow \mathbb{C}P^{kl-1}\times
\widetilde{Gr}_{k,kl}.
$$
The commutativity of the following diagram
\begin{equation}
\begin{array}{ccc}
\widetilde{Gr}_{k,kl} & \stackrel{\varphi_{k,l}}{\rightarrow} &
\mathbb{C}P^{kl-1}\times \widetilde{Gr}_{k,kl} \\
{\scriptstyle \varphi_{k,l}}\downarrow \quad \, && \qquad \quad \downarrow
{\scriptstyle \id_{\mathbb{C}P}\times \varphi_{k,l}} \\
\mathbb{C}P^{kl-1}\times \widetilde{Gr}_{k,kl} &
\stackrel{\diag_{\mathbb{C}P}\times \id_{\widetilde{Gr}}}{\rightarrow}
 & \mathbb{C}P^{kl-1}\times \mathbb{C}P^{kl-1}\times \widetilde{Gr}_{k,kl} \\
\end{array}
\end{equation}
allows us to define on the algebra $H[[x,y]]$ the structure of
$R[[z]]=\Omega^*_U(\mathbb{C}P^\infty)$-module
such that $z$ acts as the multiplication by ${\frak G}(x,y).$
Let us denote this $R[[z]]$-module by $(H[[x,y]];\: {\frak G}(x,y)).$

Let us consider $R[[z]]=\Omega^*_U(\mathbb{C}P^\infty)$ as a
Hopf algebra. Recall that $\Delta_{R[[z]]}(z)=F(z\otimes 1,1\otimes z).$

\begin{proposition}
$H[[x,y]]$ is the module coalgebra over $R[[z]],$
i. e. $R[[z]]{\mathop{\widehat{\otimes}}\limits_R}H[[x,y]]\rightarrow
H[[x,y]]$ is the homomorphism of coalgebras.
\end{proposition}
{\raggedright {\it Proof}.}\quad The proof follows from the
following commutative diagram
($(km,ln)=~1$):
\begin{equation}
\nonumber
\begin{array}{ccc}
\diagram
& {\scriptstyle \mathbb{C}P^{kl-1}\times \widetilde{Gr}_{k,kl}\times
\mathbb{C}P^{mn-1}\times \widetilde{Gr}_{m,mn}} \dlto & \\
{\scriptstyle\mathbb{C}P^{kl-1}\times \mathbb{C}P^{mn-1}\times
\widetilde{Gr}_{k,kl}\times \widetilde{Gr}_{m,mn}} \dto &&
{\scriptstyle\widetilde{Gr}_{k,kl}\times
\widetilde{Gr}_{m,mn}} \ulto_{\; \varphi_{k,l}\times \varphi_{m,n}}\dto & \\
{\scriptstyle\mathbb{C}P^{klmn-1}\times \widetilde{Gr}_{km,klmn}} &&
\quad{\scriptstyle\widetilde{Gr}_{km,klmn}}.\quad \square \llto_{\varphi_{km,ln}}
\enddiagram
\end{array}
\end{equation}

Let us consider the next commutative diagram ($(km,ln)=1$):
\begin{equation}
\begin{array}{ccc}
G\widetilde{r_{k,kl}\times Gr}_{m,mn} & \stackrel{\psi_{kl,mn}}
{\rightarrow} & \widetilde{Gr}_{km,klmn} \\
\downarrow && \downarrow \\
Gr_{k,kl}\times Gr_{m,mn} & \stackrel{\phi_{kl,mn}}{\rightarrow} &
\; Gr_{km,klmn}\; ,
\end{array}
\end{equation}
where $G\widetilde{r_{k,kl}\times Gr}_{m,mn}$ is the FBSP
over $Gr_{k,kl}\times Gr_{m,mn},$ induced by the map $\phi_{kl,mn}$
(the definition of $\phi_{kl,mn}$ was given in \cite{pap1}).
Clearly that the bundle $G\widetilde{r_{k,kl}\times Gr}_{m,mn}$
(with fiber $\mathbb{C}P^{km-1}\times \mathbb{C}P^{ln-1}$)
is (``external``) Segre's product of the canonical FBSP over $Gr_{k,kl}$
and $Gr_{m,mn}.$ By definition, put
$$
\widetilde{\; Gr\times
Gr}=\varinjlim_{(km,ln)=1}G\widetilde{r_{k,kl}\times Gr}_{m,mn}\; ,
$$
$$
\psi=\varinjlim_{(km,ln)=1}\psi_{km,ln}\colon \quad \widetilde{Gr\times Gr}
\rightarrow \widetilde{Gr}.
$$

We have the homomorphism of $R[[z]]$-modules
$$
\Psi \colon (H[[x,y]];\: {\frak G}(x,y))\rightarrow (H
{\mathop{\widehat{\otimes}}\limits_R}H[[x,y]];\: (\Delta{\frak
G})(x,y))\: ,
$$
defined by the fiber map $\psi$
(recall that $\Delta$ is the comultiplication in the Hopf algebra
$H=\Omega^*_U(Gr)$).
Clearly that the restriction $\Psi \mid_H$ coincides with $\Delta.$

Let ${\cal P}^{k-1}{\mathop{\times}\limits_X}{\cal Q}^{l-1}$
be a FBSP over a finite $CW$-complex $X$
with fiber $\mathbb{C}P^{k-1}\times \mathbb{C}P^{l-1}.$
Recall (\cite{pap1}) that if $k$ and $l$ are sufficiently large then
there exist a classifying map $f_{k,l}$
and the corresponding fiber map
\begin{equation}
\begin{array}{ccc}
{\cal P}^{k-1}{\mathop{\times}\limits_X}{\cal Q}^{l-1} &
\rightarrow & \widetilde{Gr}_{k,kl} \\
\downarrow && \downarrow \\
X & \stackrel{f_{k,l}}{\rightarrow} & Gr_{k,kl}
\end{array}
\end{equation}
which are unique up to homotopy and up to fiber homotopy respectively.
Let ${\cal P}^{km-1}{\mathop{\times}\limits_X}{\cal Q}^{ln-1},\;
(km,ln)=1$
be Segre's product of
${\cal P}^{k-1}{\mathop{\times}\limits_X}{\cal Q}^{l-1}$
with the trivial FBSP $X\times \mathbb{C}P^{m-1}\times \mathbb{C}P^{n-1}.$
Let us pass to the direct limit
$${\cal P}{\mathop{\times}\limits_X}{\cal Q}=
\varinjlim_i({\cal P}^{km_i-1}{\mathop{\times}\limits_X}
{\cal Q}^{ln_i-1}),$$ where $(km_i,ln_i)=1,$
$m_i\mid m_{i+1},\: n_i\mid n_{i+1},\; m_i\, ,n_i\rightarrow
\infty,$ as $i\rightarrow \infty.$
The stable equivalence class of FBSP (see \cite{pap1})
${\cal P}^{k-1}{\mathop{\times}\limits_X}{\cal Q}^{l-1}$
may be unique restored by the direct limit
${\cal P}{\mathop{\times}\limits_X}{\cal Q}.$
We have also a classifying map $f=\varinjlim_{(k,l)=1}f_{k,l}$
and the corresponding fiber map
\begin{equation}
\begin{array}{ccc}
{\cal P}{\mathop{\times}\limits_X}{\cal Q} &
\rightarrow & \widetilde{Gr} \\
\downarrow && \downarrow \\
X & \stackrel{f}{\rightarrow} & Gr\; .
\end{array}
\end{equation}

Let us define the category $\frak{FBSP}_f$
by the following way.
\begin{itemize}
\item[(i)] $\Ob (\frak{FBSP}_f)$ is the class of direct limits
${\cal P}{\mathop{\times}\limits_X}{\cal Q}$ of FBSP over finite
$CW$-complexes $X$ (in other words, the class of stable equivalence
classes of FBSP);
\item[(ii)] $\Mor_{\frak{FBSP}_f}({\cal P}{\mathop{\times}\limits_X}
{\cal Q},\: {\cal P}'{\mathop{\times}\limits_Y}{\cal Q}')$
is the set of fiber maps
\begin{equation}
\begin{array}{ccc}
{\cal P}{\mathop{\times}\limits_X}{\cal Q} &
\rightarrow & {\cal P}'{\mathop{\times}\limits_Y}{\cal Q}' \\
\downarrow && \downarrow \\
X & \rightarrow & Y
\end{array}
\end{equation}
such that its restrictions to any fiber
$(\cong \mathbb{C}P^\infty\times \mathbb{C}P^\infty)$ are isomorphisms.
\end{itemize}

Applying the functor of unitary cobordisms $\Omega^*_U$
to an object ${\cal P}{\mathop{\times}\limits_X}{\cal Q}
\in \Ob(\frak{FBSP}_f),$
we get the $R[[z]]$-module $(A[[x,y]];\:(f^*{\frak G})(x,y))\in
\Ob(\Omega^*_U(\frak{FBSP}_f)),$ where $A=\Omega^*_U(X)$ and $f
\colon X\rightarrow Gr$ is a classifying map for
${\cal P}{\mathop{\times}\limits_X}{\cal Q}.$
It is clear that $((\varepsilon_A\circ f^*){\frak G})(x,y)=F(x,y),$
where $\varepsilon_A \colon A\rightarrow R$ is the homomorphism,
induced by an embedding of a point $\pt \hookrightarrow X.$
In other words, for any object in the category $\Omega^*_U(\frak{FBSP_f})$
there exists the canonical morphism $(A[[x,y]];\:(f^*{\frak G})(x,y))
\rightarrow (R[[x,y]];\: F(x,y)).$

Hence there exist the initial object $(H[[x,y]];\: {\frak G}(x,y))$
and the final object $(R[[x,y]];\: F(x,y))$ in the category
$\Omega^*_U(\frak{FBSP}).$

Let's consider a pair $(A[[x,y]];\:(f^*{\frak G})(x,y)),\;
(B[[x,y]];\:(g^*{\frak G})(x,y))\in \Ob(\Omega^*_U(\frak{FBSP_f})),$
where $(B[[x,y]];\:(g^*{\frak G})(x,y))=
\Omega^*_U({\cal P}'{\mathop{\times}\limits_Y}{\cal Q}').$
Let's define their ``tensor product`` as the object
$((A{\mathop{\otimes}\limits_R}B)[[x,y]];\:(((f^*\otimes g^*)\circ
\Delta){\frak G})(x,y))\in \Ob(\Omega^*_U(\frak{FBSP_f}))$
(recall that $\Delta \colon H\rightarrow
H{\mathop{\widehat{\otimes}}\limits_R}H$ is the comultiplication
in the Hopf algebra $H$).

\begin{proposition}
The category $\Omega^*_U(\frak{FBSP_f})$ is the tensor category
with the just defined tensor product and the unit
$\; \I =(R[[x,y]];\: F(x,y)).$
\end{proposition}
{\raggedright {\it Proof}.}\quad The proof is trivial.
For example, the associativity axiom follows from the
identity $(((\Delta \otimes \id_H)\circ \Delta){\frak G})(x,y)=
(((\id_H\otimes \Delta)\circ \Delta){\frak G})(x,y)$
which follows from the next commutative diagram ($(kmt,lnu)=1$):
\begin{equation}
\begin{array}{ccc}
G\widetilde{r_{km,klmn}\times Gr}_{t,tu} &
\rightarrow & \widetilde{Gr}_{kmt,klmntu} \\
\uparrow && \uparrow \\
Gr_{k,kl}\widetilde{\times Gr_{m,mn}}
\times Gr_{t,tu} &
\rightarrow & G\widetilde{r_{k,kl}\times G}r_{mt,mntu},
\end{array}
\end{equation}
where $Gr_{k,kl}\widetilde{\times Gr_{m,mn}}
\times Gr_{t,tu}$
is external Segre's product of the canonical FBSP over
$Gr_{k,kl},\; Gr_{m,mn}$ and $Gr_{t,tu}$ (it
is the bundle over $Gr_{k,kl}\times Gr_{m,mn}\times Gr_{t,tu}$ with
fiber $\mathbb{C}P^{kmt-1}\times \mathbb{C}P^{lnu-1}$).\:$\square$

\smallskip

Note that there exist the canonical homomorphisms $p_1,\; p_2$:
$$
\diagram
& \scriptstyle{((A{\mathop{\otimes}\limits_R}B)[[x,y]];\:
(((f^*\otimes g^*)\circ
\Delta){\frak G})(x,y))} \dlto_{\scriptscriptstyle{p_1}}
\drto^{\scriptscriptstyle{p_2}} &\\
\scriptstyle{(A[[x,y]];\:(f^*{\frak G})(x,y))} &&
\scriptstyle{(B[[x,y]];\:(g^*{\frak
G})(x,y))}\; ,\\
\enddiagram
$$
such that $p_1\mid_{A{\mathop{\otimes}\limits_R}B}=\id_A\otimes
\varepsilon_B,\quad
p_2\mid_{A{\mathop{\otimes}\limits_R}B}=\varepsilon_A
\otimes \id_B.$


\begin{thebibliography}{99}
\bibitem{pap1}
{\sc A. V. Ershov}
Floating bundles and their applications.---
arXiv:math.AT/0102054
\end{thebibliography}
\end{document}